\documentclass[12pt,a4paper,leqno]{amsart}

\usepackage{anysize}
\usepackage[utf8]{inputenc}
\usepackage{amsmath}
\usepackage{amssymb}
\usepackage{array}

\marginsize{3cm}{3cm}{3.5cm}{3.5cm}

\newtheorem{theorem}{Theorem}[section]

\newtheorem{remark}[theorem]{Remark}

\begin{document}

\title{Self-replication and Borwein-like algorithms}
\author{Jesús Guillera}
\address{Department of Mathematics, University of Zaragoza, 50009 Zaragoza, SPAIN}
\email{jguillera@gmail.com}
\dedicatory{To the memory of Jonathan Borwein}
\date{}

\maketitle

\begin{abstract}
Using a self-replicating method, we generalize with a free parameter some Borwein algorithms for the number $\pi$. This generalization includes values of the Gamma function like $\Gamma(1/3)$, $\Gamma(1/4)$ and of course $\Gamma(1/2)=\sqrt{\pi}$. In addition, we give new rapid algorithms for the perimeter of an ellipse.
\end{abstract}

\section{Introduction}

In the mid '80s the Borwein brothers observed that the modular equations given by Ramanujan could be connected to the number $\pi$ via the elliptic alpha function, to construct extraordinarily rapid algorithms to calculate $\pi$ \cite[p. 170]{Bo}. The point of view of the Borweins is explained in \cite{Bo2}, and their proofs require a good knowledge of the elliptic modular functions and forms. In \cite{Gui2} we presented a different and simple strategy and we used it further in \cite{Coo-Gui-St-Zu} applying it to new self-replicating identities of non-hypergeometric type. It is precisely in \cite{Coo-Gui-St-Zu} where we named the technique as {\it self-replication method}, and its property as {\it Hidden Modularity}, which makes knowledge of modularity (although interesting) unnecessary to prove the algorithms. In this paper we generalize, with a free parameter, some Borwein algorithms for the number $\pi$, to  include some values of the Gamma function, like $\Gamma(1/3)$, $\Gamma(1/4)$ and of course $\Gamma(1/2)=\sqrt{\pi}$. In addition, we give new rapid algorithms for the perimeter of an ellipse.

\section{Algebraic transformations}

In this section we use known algebraic transformations to obtain self-replication identities. 

\subsection{Landen's transformation}

Landen (1719-1790) proved the following remarkable formula \cite[p. 591]{Almkvist-Berndt}:
\[
\sum_{k=0}^{\infty} \frac{\left(a\right)_k\left(b\right)_k}{(2b)_k \left(1\right)_k} \left( \frac{4t}{(1+t)^2}\right)^k 
= (1+t)^{2a} \sum_{k=0}^{\infty} \frac{\left(a)_k(a-b+\frac12\right)_k}{\left(b+\frac12 \right)_k \left(1\right)_k} t^{2k}.
\]
Here, we only will need the particular case $a=b=1/2$ written in the following form \cite[Thm 1.2 (b)]{Bo}:
\begin{equation}\label{quad-trans}
\sum_{k=0}^{\infty} \frac{\left(\frac12\right)_k^2}{(1)_k^2}x^{2k}=(1+t) \sum_{k=0}^{\infty} \frac{\left(\frac12\right)_k^2}{(1)_k^2}t^{2k}, \qquad t = \frac{1-\sqrt{1-x^2}}{1+\sqrt{1-x^2}}.
\end{equation}
Observe that this transformation is quadratic because $t \sim 4^{-1} \cdot x^2$. In \cite{Gui} we explain how to re-prove known algebraic transformations in an elementary way. Then, if we apply the differential operator
\[
a+\frac{b}{2} \, \vartheta_x = a+b \, \frac{1+t}{1-t} \, \vartheta_t,
\]
to both sides of (\ref{quad-trans}), we have
\[
\left( a+\frac{b}{2} \, \vartheta_x \right) \left\{ \sum_{k=0}^{\infty} \frac{\left(\frac12\right)_k^2}{(1)_k^2}x^{2k} \right\} =
\left( a+b \, \frac{1+t}{1-t} \, \vartheta_t \right) \left\{ (1+t) \sum_{k=0}^{\infty} \frac{\left(\frac12\right)_k^2}{(1)_k^2}t^{2k} \right\},
\]
and we arrive at the self-replicating identity:
\begin{equation}\label{comp-trans}
\sum_{k=0}^{\infty} \frac{\left(\frac12\right)_k^2}{(1)_k^2} (a+bk) x^{2k}=
\sum_{k=0}^{\infty} \frac{\left(\frac12\right)_k^2}{(1)_k^2} 
\left( a(1+t)+b\frac{t(1+t)}{1-t} + 2b\frac{(1+t)^2}{1-t} k \right) t^{2k}.
\end{equation}
Writing 
\[ 
\alpha= a(1+t)+b\frac{t(1+t)}{1-t}, \qquad \beta=2b\frac{(1+t)^2}{1-t},
\]
we see the self-replication with clarity. In Section \ref{Sect-algos} we will use the fact that the product of a power of the identity (\ref{quad-trans}) by identity (\ref{comp-trans}) also has a self-replication property.

\subsection{Cubic transformation}

The following known algebraic transformation
\begin{equation}\label{cubic-trans}
\sum_{k=0}^{\infty} \frac{\left(\frac13\right)_k\left(\frac23\right)_k}{(1)_k^2}x^{3k}=(1+2t) \sum_{k=0}^{\infty} \frac{\left(\frac13\right)_k\left(\frac23\right)_k}{(1)_k^2}t^{3k}, \qquad t=\frac{1-\sqrt[3]{1-x^3}}{1+2\sqrt[3]{1-x^3}},
\end{equation}
\cite[Chapter 4]{Cooper-book}, is of degree $3$ (cubic). Applying to it the operator
\[
a+\frac{b}{3}\vartheta_x = a + b \frac{(1+2t)(1-t^3)}{(1-t)^3} \vartheta_t,
\]
we arrive at the identity
\begin{multline}\label{comp-cubic}
\sum_{k=0}^{\infty} \frac{\left(\frac13\right)_k\left(\frac23\right)_k}{(1)_k^2} (a+bk) x^{3k} \\ =
\sum_{k=0}^{\infty} \frac{\left(\frac13\right)_k\left(\frac23\right)_k}{(1)_k^2} 
\left( a(1+2t)+2b\frac{t(1+2t) (1-t^3)}{(1-t)^3} + 3b\frac{(1-t^3)(1+2t)^2}{(1-t)^3} k \right) t^{3k}.
\end{multline}
Writing
\[
\alpha=a(1+2t)+2b\frac{t(1+2t) (1-t^3)}{(1-t)^3}, \qquad \beta=3b\frac{(1-t^3)(1+2t)^2}{(1-t)^3},
\]
we see that it is a self-replicating identity.

\subsection{Quartic transformation} 

Finally, consider the following known quartic algebraic transformation \cite[Chapter 4]{Cooper-book}:
\begin{equation}\label{quartic-trans}
\sum_{k=0}^{\infty} \frac{\left(\frac12\right)_k^2}{(1)_k^2}x^{4k}=(1+t)^2 \sum_{k=0}^{\infty} \frac{\left(\frac12\right)_k^2}{(1)_k^2}t^{4k}, \qquad t = \frac{1-\sqrt[4]{1-x^4}}{1+\sqrt[4]{1-x^4}}.
\end{equation}
Applying to both sides the operator
\[
a+\frac{b}{4} \vartheta_x = a+b\frac{(1+t^2)(1+t)}{(1-t)^3}\vartheta_t,
\]
we have the identity
\begin{multline}\label{comp-quartic}
\sum_{k=0}^{\infty} \frac{\left(\frac12\right)_k^2}{(1)_k^2} (a+bk) x^{4k} = \\
\sum_{k=0}^{\infty} \frac{\left(\frac12\right)_k^2}{(1)_k^2} 
\left( a(1+t)^2+2b\frac{t(1+t^2)(1+t)^2}{(1-t)^3} + 4 b\frac{(1+t^2) (1+t)^3}{(1-t)^3} k \right) t^{4k}.
\end{multline}
In this case writing
\[
\alpha = a(1+t)^2+2b\frac{t(1+t^2)(1+t)^2}{(1-t)^3}, \qquad \beta = 4 b\frac{(1+t^2) (1+t)^3}{(1-t)^3},
\]
we see that it is a self-replicating identity.

\section{A generalization of Borweins' algorithms for $\pi$}\label{Sect-algos} 

In this section we use self-replication formulas and couples of series of Ramanujan-type to derive a kind of rapid algorithms of Borwein-type. A couple  of such series is the following known one:

\begin{align*}
\sum_{k=0}^{\infty} \frac{(s)_k(1-s)_k}{(1)_k^2} \frac{1}{2^k} &= 
\frac{\sqrt{\pi}}{\Gamma \left(1-\frac{s}{2} \right) \Gamma \left(\frac12+\frac{s}{2}\right)}, 
\\
\sum_{k=0}^{\infty} \frac{(s)_k(1-s)_k}{(1)_k^2} \frac{k}{2^k} &= 
\frac{s \, \sqrt{\pi}}{\Gamma \left(1+\frac{s}{2} \right) \Gamma \left(\frac12-\frac{s}{2}\right)},
\end{align*}
which can be proved automatically by the Wilf-Zeilberger (WZ)-method due to the presence of the free parameter $s$. We will use this couple of identities to get initial values for our algorithms, but observe that this is only possible for $s=1/2, \, 1/4, \, 1/3, \, 1/6$, because we need to use algebraic transformations, to derive the recurrences of the algorithms, and they do not exist for other values of $s$. In \cite[Sect. 4]{Gui2} we see other couples of series which can be used to get the initial values and the result to which the quadratic and quartic algorithms tend. To deduce the algorithms we will use the method explained in \cite{Gui2}.

\subsection{A quadratic algorithm}

Taking $s=1/2$, we get
\[
\sum_{k=0}^{\infty} \frac{\left(\frac12\right)_k^2}{(1)_k^2} \frac{1}{2^k} = \frac{\sqrt{\pi}}{\Gamma^2 \left( \frac34 \right)}, \qquad 
\sum_{k=0}^{\infty} \frac{\left(\frac12\right)_k^2}{(1)_k^2} \frac{k}{2^k} = \frac{\Gamma^2 \! \left( \frac34 \right)}{{\pi^{3/2}}}.
\]
Hence, we have the following formula:
\begin{equation}\label{initial-gamma34}
\left( \sum_{k=0}^{\infty} \frac{\left(\frac12\right)_k^2}{(1)_k^2} \frac{1}{2^k} \right)^{\! \! w}
\sum_{k=0}^{\infty} \frac{\left(\frac12\right)_k^2}{(1)_k^2} \frac{k}{2^k} =  \frac{1}{\Gamma\left(\frac34\right)^{2w-2}\pi^{\frac32-\frac{w}{2}}}.
\end{equation}
From (\ref{quad-trans}) and (\ref{comp-trans}), letting
\[
A_n = \left( \sum_{k=0}^{\infty} \frac{\left(\frac12\right)_k^2}{(1)_k^2}  d_n^{2k} \right)^{\! \! w} \sum_{k=0}^{\infty} \frac{\left(\frac12\right)_k^2}{(1)_k^2} (a_n+b_nk) d_n^{2k},
\]
observing that $\lim a_n = \lim A_n = A_0$ (see \cite{Gui2} for the explanation), getting the initial values and the limit from (\ref{initial-gamma34}), and defining $c_n=b_n/(1-d_n^2)$, we arrive at
\begin{align*}
& d_0=\frac{1}{\sqrt{2}}, \quad c_0=2, \quad a_0=0, \\
& d_{n+1}=\frac{1-\sqrt{1-d_n^2}}{1+\sqrt{1-d_n^2}}, \quad c_{n+1}=2c_n(1+d_{n+1})^{w-1}, \\
& a_{n+1}=a_n(1+d_{n+1})^{w+1}+\frac12 c_{n+1}d_{n+1}(1-d_{n+1}), \\
& a_n(w) \to \frac{1}{\Gamma\left(\frac34\right)^{2w-2}{\pi^{\frac32-\frac{w}{2}}}},
\end{align*}
which is a quadratic algorithm. Some examples are
\[
a_n(1) \to \frac{1}{\pi}, \quad a_n(3) \to \frac{1}{\Gamma^4 \! \left(\frac34\right)}, \quad a_n(\frac13) \to \left( \frac{\sqrt 2}{\Gamma\left(\frac14\right)} \right)^{\! \! \frac43},
\]
which are algebraic algorithms for $\pi$, $\Gamma(3/4)$ and $\Gamma(1/4)$ respectively. We see that the case $w=1$ with the substitution $a_n=r_n- 2^n \, d_n^2$, is a Borweins' quadratic algorithm for $\pi$. A related algorithm for $\Gamma(1/4)$ is given in \cite[p. 137]{Ba-Bo}. In \cite{Gui} we deduced the Borweins' quadratic algorithm for $\pi$ from the Gauss-Salamin-Brent algorithm, and it is easy to prove that in fact both algorithms are equivalent.

\subsection{A cubic algorithm}

In a similar way, from the (\ref{cubic-trans}) and (\ref{comp-cubic}), defining
\[
A_n = \left( \sum_{k=0}^{\infty} \frac{\left(\frac13\right)_k\left(\frac23\right)_k}{(1)_k^2}  d_n^{3k} \right)^{\! \! w} \sum_{k=0}^{\infty} \frac{\left(\frac13\right)_k \left(\frac23\right)_k}{(1)_k^2} (a_n+b_nk) d_n^{3k},
\]
observing that $\lim a_n = \lim A_n = A_0$, letting $c_n=b_n/(1-d_n^3)$, and taking the initial values from
\begin{equation}\label{initial-gamma23}
\left( \sum_{k=0}^{\infty} \frac{\left(\frac13\right)_k\left(\frac23\right)_k}{(1)_k^2} \frac{1}{2^k} \right)^{\! \! w}
\sum_{k=0}^{\infty} \frac{\left(\frac13\right)_k\left(\frac23\right)_k}{(1)_k^2} \frac{k}{2^k} = 
\frac{3^{1-\frac{w}{2}}\cdot 2^{\frac{2w}{3}-\frac53}}{\pi^{2-w} \cdot \Gamma \left(\frac23\right)^{3w-3}},
\end{equation}
we arrive at the following algorithm
\begin{align*}
& d_0=\frac{1}{\sqrt[3]{2}}, \quad c_0=2, \quad a_0=0, \\
& d_{n+1}=\frac{1-\sqrt[3]{1-d_n^3}}{1+2\sqrt[3]{1-d_n^3}}, \quad c_{n+1}=3c_n(1+2d_{n+1})^{w-1}, \\
& a_{n+1}=a_n(1+2d_{n+1})^{w+1}+\frac23 c_{n+1}d_{n+1} \frac{1-d_{n+1}^3}{1+2d_{n+1}}, \\ 
& a_n(w) \to \frac{3^{1-\frac{w}{2}}\cdot 2^{\frac{2w}{3}-\frac53}}{\pi^{2-w} \cdot \Gamma \left(\frac23\right)^{3w-3}},
\end{align*}
which generalizes a cubic algorithm for $\pi$ due to the Borwein brothers. Some examples are
\[
a_n(1) \to \frac{\sqrt 3}{2 \pi}, \quad a_n(2) \to \frac{1}{\sqrt[3]{2} \cdot \Gamma^3 \! \left(\frac23\right)}, \quad a_n(\frac12) \to \left( \frac{2}{\sqrt{3} \Gamma\left(\frac13\right)} \right)^{\! \! \frac32},
\]
which are algebraic algorithms for $\pi$, $\Gamma(2/3)$ and $\Gamma(1/3)$ respectively. 

\subsection{A quartic algorithm}

From (\ref{quartic-trans}) and (\ref{comp-quartic}), defining
\[
A_n = \left( \sum_{k=0}^{\infty} \frac{\left(\frac12\right)_k^2}{(1)_k^2}  d_n^{4k} \right)^{\! \! w} \sum_{k=0}^{\infty} \frac{\left(\frac12\right)_k^2}{(1)_k^2} (a_n+b_nk) d_n^{4k}, 
\]
letting $c_n=b_n/(1-d_n^4)$ and getting the initial values from (\ref{initial-gamma34}), we see that taking
\[
d_0=\frac{1}{\sqrt[4]{2}}, \quad c_0=2, \quad a_0=0,
\]
and the recurrences

\begin{align}\label{recu-quartic}
& d_{n+1}=\frac{1-\sqrt[4]{1-d_n^4}}{1+\sqrt[4]{1-d_n^4}}, \quad c_{n+1}=4c_n(1+d_{n+1})^{2w-2}, \\
& a_{n+1}=a_n(1+d_{n+1})^{2w+2}+\frac12 c_{n+1} \frac{d_{n+1}}{1+d_{n+1}} (1-d_{n+1}^4), 
\end{align}
we have
\[
a_n(w) \to \frac{1}{\Gamma\left(\frac34\right)^{2w-2}{\pi^{\frac32-\frac{w}{2}}}}.
\]
Some examples are
\[
a_n(1) \to \frac{1}{\pi}, \quad a_n(3) \to \frac{1}{\Gamma^4 \! \left(\frac34\right)}, \quad a_n(\frac13) \to \left( \frac{\sqrt 2}{\Gamma\left(\frac14\right)} \right)^{\! \! \frac43}.
\]
In \cite[Chapter 31]{Arndt} there are many algorithms for $\pi$, $\Gamma(1/4)$, $\Gamma(1/8)$ and other constants related to the evaluation of the elliptic integral $K$ at the singular values. 

\section{The perimeter of an ellipse}

The perimeter $P(a,b)$ of an ellipse of semi-axis $a$ and $b$ with $a \geq b$, is given by the elliptic integral 
\[
P(a,b)=4a\int_{0}^{\frac{\pi}{2}} \sqrt{1-x^2 \cos^2 \varphi} \, \,  d \varphi, \qquad x=\sqrt{1-\frac{b^2}{a^2}},
\]
where $x$ is the eccentricity of the ellipse. The development in series of powers of $x$ leads to the following hypergeometric formula \cite[p. 38-39]{ramanujan}:
\begin{equation}\label{P-ellipse}
P(a,b)=\frac{2\pi b^2}{a} \sum_{k=0}^{\infty} \frac{\left(\frac12\right)_k^2}{(1)_k^2} (1+2k) \left( 1-\frac{b^2}{a^2} \right)^{\! \! k},
\end{equation}
which is known to be related to the AGM (Arithmetic Geometric Mean) of Gauss. For details and an historical account see \cite{Almkvist-Berndt}. The authors of that paper comment that an Ivory's letter of 1796 unmistakenly pointed that he knew it. The following algorithms come from our method of self-replication instead of using the AGM. 

\subsection{A quadratic algorithm for the perimeter}

From the hypergeometric formula (\ref{P-ellipse}) for $P(a,b)$ and the self-replicating transformation (\ref{comp-trans}), we can arrive at the following quadratic algorithm:
\begin{align*}
& d_0=\sqrt{1-\frac{b^2}{a^2}}, \qquad a_0=1, \qquad c_0=\frac{2a^2}{b^2}, \\
& d_{n+1} = \frac{1-\sqrt{1-d_n^2}}{1+\sqrt{1-d_n^2}}, \qquad c_{n+1}=\frac{2c_n}{1+d_{n+1}}, \\
& a_{n+1} = a_n(1+d_{n+1}) + \frac12 c_{n+1} d_{n+1}(1-d_{n+1}), \qquad \frac{2\pi b^2}{a} \, a_n \to P(a,b)
\end{align*}
In \cite[eq. 27]{Almkvist-Berndt} there is a related quadratic algorithm for $P(a,b)$ based on the AGM. 

\subsection{A quartic algorithm for the perimeter}

From the quartic recurrences (\ref{recu-quartic}) with $w=0$:
\begin{align*}
& d_{n+1}=\frac{1-\sqrt[4]{1-d_n^4}}{1+\sqrt[4]{1-d_n^4}}, \quad c_{n+1}=4 \frac{c_n}{(1+d_{n+1})^2}, \\
& a_{n+1}=a_n(1+d_{n+1})^2+\frac12 c_{n+1} \frac{d_{n+1}}{1+d_{n+1}} (1-d_{n+1}^4), 
\end{align*}
and taking as initial values 
\[
d_0=\sqrt[4]{1-\frac{b^2}{a^2}}, \qquad c_0=\frac{b_0}{1-d_0^4}=\frac{2a^2}{b^2}, \qquad a_0=1,
\]
we see that 
\[
\frac{2\pi b^2}{a} a_n \to P(a,b),
\]
quartically, that is multiplying by $4$ the number of correct digits in each iteration.

\begin{remark} \rm
The papers \cite{Almkvist-Berndt} and \cite{Bo2} are reprinted in \cite{Ba-Bo2}, and the paper \cite{ramanujan} is reprinted in \cite{BeBo}. The books \cite{BeBo} and \cite{Ba-Bo2} are two collections of very interesting papers related to the number $\pi$. The book \cite{Arndt} will be useful for the computationalist, whether a working programmer or anyone interested in methods of computation. The recent nice book \cite{Cooper-book} by Shaun Cooper deals with Ramanujan's theta functions, and in the Chapter $14$ there are applications to series and algorithms for $1/\pi$.
\end{remark}

\end{document}